\documentclass[12pt,reqno]{amsart}
\usepackage{latexsym,amsmath,amsthm,amssymb,amsfonts}
\usepackage{setspace}
\usepackage{color}
\usepackage{tikz}

\numberwithin{equation}{section}
\newtheorem{theorem}{Theorem}[section]
\newtheorem{lem}[theorem]{Lemma}
\newtheorem{thm}[theorem]{Theorem}
\newtheorem{pro}[theorem]{Proposition}
\newtheorem{cor}[theorem]{Corollary}

\def\S2{\mathbb{S}^2}
\def\s{\,\,\,\,}

\def\dint{\displaystyle{\int}}
\def\mv{1.7ex}
\def\endproof{$\hfill\Box$\\}
\def\R{\mathbb{R}}
\def\C{\mathbb{C}}

\title{\bf Metrics  On $S^2$ With  bounded $\|K_g\|_{L^1\log L^1}$ and small $\|K_g-1\|_{L^1}$}
\author{Yuxiang Li, Hongyan Tang}

\address{\newline
Yuxiang Li:
 Department of Mathematical Sciences, Tsinghua University, Beijing 100084, P.R. China.
{\tt Email:yxli@math.tsinghua.edu.cn}
\newline
\newline
Hongyan Tang:
Department of Mathematical Sciences, Tsinghua University, Beijing 100084, P.R. China.
{\tt Email:hytang@math.tsinghua.edu.cn}}

\date{}
\begin{document}
\maketitle

\begin{abstract}
In this short paper, we will study the convergence of 
a metric sequence $g_k$
on $S^2$ with bounded $\int|K_{g_k}|\log(1+|K_{g_k}|)d\mu_{g_k}$ and
small $\int|K(g_k)-1\|d\mu_{g_k}$. We will show that 
 such a sequence is precompact.

\end{abstract}

\section{Introduction}
Let $(M,g)$ be a closed Riemannian manifold with $\dim M\geq 3$,
and $g_k=u_k^\frac{4}{n-2}g$ be a metric.
Gursky \cite{Gursky} proved   that 
a subsequence of $g_k$ converges in $C^{0,\alpha}$ whenever 
$vol(g_k)=1$ and $\|K_{g_k}\|_{L^p(M,g_k)}<C$, where $p>\frac{n}{2}$.  
In \cite{Li-Zhou,Li-Zhou-Wei}, we showed that Gursky's result can be deduced
from the fact that $S^{n-1}\times R$ does not carry a flat metric
which is conformal to $dt^2+g_{\mathbb{S}^{n-1}}$.
However, such a fact does not hold for the case of $n=2$, and
therefore Gursky's
result is not true for 2 dimensional case.

For example, put $Q_k= S^1\times (-k\pi^2+2\pi,k\pi^2-2\pi)$,
and 
\begin{equation}\label{metric}
g_k=\left(\frac{1}{2k\pi\cos(\frac{t}{2k\pi})}\right)^2(dt^2+d\theta^2).
\end{equation}
We have
$$
K_{g_k}=-\left(2k\pi\cos(\frac{t}{2k\pi})\right)^2\Delta
\frac{1}{2}\log\left(\frac{1}{2k\pi\cos(\frac{t}{2k\pi})}\right)^2=-1,
$$ 
and 
$$
\mu(g_k,Q_k)=2\pi\int_{-k\pi^2+2\pi}^{k\pi^2-2\pi}\left(\frac{1}{2k\pi\cos(\frac{t}{2k\pi})}\right)^2dt
=2\frac{\cos\frac{1}{k}}{k\sin\frac{1}{k}}.
$$
Noting that on $[-k\pi^2+2\pi,-k\pi^2+4\pi]$, $g_k\sim \frac{dt^2+d\theta^2}{s^2}$, where $s=t-(-k\pi^2+2\pi)$,
for any $p>1$, we can easily extend $(Q_k,g_k)$ to a sphere with a smooth 
metric $h_k$, such that 
$$
\|K_{h_k}\|_{L^p(S^2,h_k)}+\mu(h_k,S^2)<C.
$$  
Since the conformal class of $S^2$ is unique, we may assume 
$h_k$ is conformal to $g_{\S2}$. Obviously, $h_k$ can not converge
in $C^{0,\alpha}$. In fact the diameter of $(S^2,h_k)$ tends to
infinity.

\begin{center}
\begin{tikzpicture}

\draw (0,0) arc (50: 130: 3);
\draw (-3.8,1.6) arc (-130: -50: 3);

\draw[dotted, line width=1.3pt] (0,0) arc (-130: 130: 1.02);
\draw[dotted, line width=1.3pt] (-3.8,1.6) arc (50: 310: 1.02);

\end{tikzpicture}

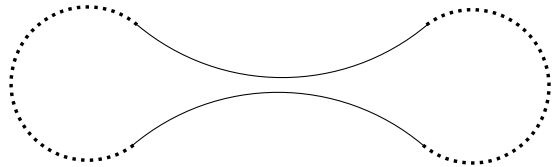
{\bf Fig 1. $(S,h_k)$}
\end{center}

Furthermore, we can define a $W^{2,p}$-metric on 
$Q'_k=S^1\times(-k\pi^2+2\pi-k^2,k\pi^2-2\pi+k^2)$ as follows:
$$
g_k'=\left\{\begin{array}{lll}
\left(\frac{1}{2k\pi\cos(\frac{t+k^2}{2k\pi})}\right)^2(dt^2+d\theta^2)&&
t<-k^2\\[\mv]
\left(\frac{1}{2k\pi}\right)^2(dt^2+d\theta^2)&&t\in[-k^2,k^2]\\[\mv]
\left(\frac{1}{2k\pi\cos(\frac{t-k^2}{2k\pi})}\right)^2(dt^2+d\theta^2)&&t>k^2.
\end{array}\right.
$$
Then we have $K_{g_k'}=-1$ or $0$, and
$$
\mu(g_k',Q_k')=2\frac{\cos\frac{1}{k}}{k\sin\frac{1}{k}}+\frac{1}{2\pi^2}.
$$
We can also extend $(Q_k',g_k')$ to a sphere with a $W^{2,p}-$ 
metric $h_k'$, where
$$
\|K_{h_k'}\|_{L^p(S^2,h_k')}+\mu(h_k',S^2)<C.
$$ 
We can approximate $h_k'$ in $W^{2,p}$ by smooth metrics. Then we can
find a smooth metric $\hat{h}_k$, which
does not converge in $C^{0,\alpha}$, and 
whose eara does not convege to the eara of the Gromov-Hausdorff
limit.

\begin{center}
\begin{tikzpicture}

\draw (0,0) arc (50: 90: 3);
\draw (-1.9,0.9) arc (-90: -50: 3);
\draw[dashed] (-1.9,0.9)--(-4,0.9);
\draw[dashed] (-1.9,0.7)--(-4,0.7);

\draw (-4,0.7) arc (90: 130: 3);
\draw (-5.9,1.61) arc (-130: -90: 3);

\draw[dotted, line width=1.3pt] (0,0) arc (-130: 130: 1.02);
\draw[dotted, line width=1.3pt] (-5.9,1.61) arc (50: 310: 1.02);

\end{tikzpicture}

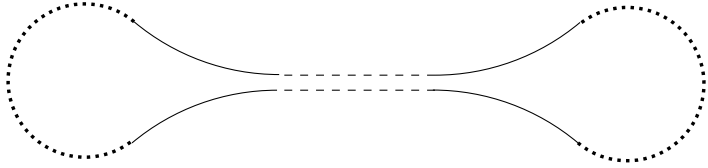
{\bf Fig 2. $(S,h_k')$}
\end{center}
For more examples and  results, one can refers to \cite{Chen}.

\vspace{1ex}
In this short paper, we will prove that the convergence holds when
$\|K_g-1\|_{L^p}$ small. In fact, we will consider a more general case.

Let
$$
\mathcal{M}_\Lambda=
\{g=e^{2u}g_{\mathbb{S}^2}:
\int_{S^2}|K_g|\log(1+|K_g|)d\mu_g<\Lambda\}.
$$
Here
$$
\S2=\left\{(x^1,x^2,x^3)\in\R^3:\sum_{i=1}^3(x^i)^2=1\right\}.
$$
The main result of this paper is the following:

\begin{thm} \label{main} 
There exist positve constants $\tau$ and $\lambda$, such that 
if $g=e^{2u}g_{\S2}\in\mathcal{M}_\Lambda$ and
$$
\int_{S^2}|K_g-1|d\mu_g<\tau,
$$
then there exists a Mobius transformation 
$\sigma$, such that $\sigma^*(g)=e^{2u'}$
with  $|u'|<\lambda$.
\end{thm}

We set
$$
\mathcal{M}^p(\Lambda, \lambda)=
\{u: g_u=e^{2u}g_{\mathbb{S}^2},\s \mu(g_u)\leq \Lambda,\s
\int_{S^2}|K_{g_u}-1|^pd\mu_{g_u}<\lambda \}.
$$
Noting that $|K|\log(1+|K|)<C(p)(1+|K|^p)$ when $p>1$, we have
\begin{thm}
Let $p>1$. 
Then there exists $\lambda_0=\lambda_0(p,\Lambda)$, 
such that  $\mathcal{M}^p(\Lambda,\lambda)$
is weakly compact in $W^{2,p}$ up to M\"obius transformations, whenever $\lambda<\lambda_0$,
\end{thm}

\section{Some previous results }
The classification of solutions of equation 
$$
-\Delta u=e^{2u}
$$
on $\R^2$ with $\int_{\R^2}e^{2u}<+\infty$
has been solved by Wenxiong Chen and Congming Li \cite{Chen-Li}. They proved the following:
\begin{thm}\label{Chen-Li}
Let $-\Delta u=e^{2u}$. If $\int_{\R^2}e^{2u}<+\infty$, then
$$
u=-\log(1+\frac{1}{4}|x-x_0|^2),\s and\s \int_{\R^2}e^{2u}=4\pi.
$$
\end{thm}

A special case of the above theorem can be deduced from Ding's lemma (\cite{Ding}, c.f. \cite{Chen-Li} Lemma 1.1)
\begin{lem}\label{Ding}
Let $-\Delta u=e^{2u}$. If $\int_{\R^2}e^{2u}\leq 4\pi$, then
$$
u=-\log(1+\frac{1}{4}|x-x_0|^2),\s and\s \int_{\R^2}e^{2u}=4\pi.
$$
\end{lem}

To get a Harnack estimate on the neck domain, we  will use 
the following result  (for the proof, one can refer to \cite{Kuwert-Li}):

\begin{lem}\label{Lq}
Let $(\Sigma,g)$ be
a closed Riemann surface.
If $u$ solves the equation
$$-\Delta u =f,$$
then for any $q\in (0,2)$ $r>0$ and $x\in \Sigma$, 
$$
r^{2-q}\int_{B_r(x)}|\nabla_gu|^q\leq C(\Sigma,g,q)\|f\|_{L^1(\Sigma)}.
$$
\end{lem}

We set
$$L^1{\log L^1}(\Omega)=\{f\in L^1(\Omega):\int_\Omega|f|\log(1+|f|)<+\infty\},$$
and
$$\|f\|_{L^1{\log L^1}(\Omega)}=\int_\Omega|f|\log(1+|f|)dx.$$
We have
\begin{lem} Let $I(f)=\int_{\R^2}f(y)\log|x-y|dy$. Then
$$\|I(f)\|_{L^\infty(D_R)}
\leq C(R)(1+\|f\|_{L^1{\log L^1}(\R^2)}).$$
\end{lem}

\proof Let $x\in D_R$. We have
$$|I(f)|(x)\leq\int_{D_\frac{1}{2}(x)}|f(y)\log |x-y||dy+\int_{D_{2R}\setminus
    D_\frac{1}{2}(x)}|f(y)\log |x-y||dy.$$
Since
$$\begin{array}{lll}
  \dint_{D_\frac{1}{2}(x)}|f(y)\log |x-y||dy&=&\dint_{\{y\in D_\frac{1}{2}(x),\,\,1+|f(y)|\geq \frac{1}{\sqrt{|x-
  y|}}\}}
   \cdots+\int_{\{y\in D_\frac{1}{2}(x),\,\,1+|f(y)|< \frac{1}{\sqrt{|x-y|}}\}}\cdots\\[\mv]
   &\leq&C\dint_{D_\frac{1}{2}(x)}|f(y)|\log(1+|f|)dy\\[\mv]
   &&+\dint_{D_\frac{1}{2}(x)}|\log|x-y||\left(\frac{1}{|x-y|^\frac{1}{2}}
   -1\right)dy,
  \end{array}$$
and
$$\int_{D_{2R}\setminus D(x)}|f(y)\log |x-y||dy\leq\log 2R
\int_{\R^2}|f(y)|dy\leq C(R)(1+\|f\|_{L^1\log L^1(D_{2R})}),$$
we get
$$\|I(f)\|_{L^\infty(D_R)}
\leq C(R)(1+\|f\|_{L^1{\log L^1}(\R^2)}).$$

\endproof

\begin{cor}\label{Linfinity2} Let $f\in L^1{\log L^1}(D)$
and $u$ solve the equation
\begin{equation*}
-\Delta u=f,\s u|_{\partial D}=0.
\end{equation*}
Then
\begin{equation}\label{Hardy2.1}
\|u\|_{L^\infty(D)}
\leq C(1+\|f\|_{L^1{\log L^1}(D)})
\end{equation}
and
\begin{equation}\label{Hardy2.2}
\|\nabla u\|_{L^2(D)}
\leq C(1+\|f\|_{L^1{\log L^1}(D)}).
\end{equation}
\end{cor}

\proof Let
$\tilde{f}=\left\{\begin{array}{ll}
f&x\in D\\
  0&x\notin D
\end{array}\right.$ and $v=I(\tilde{f})$.
Then $v-u$ is harmonic. By Maximum Principle, we get
 \eqref{Hardy2.1}.

If we let $f_k=\max\{\min\{f,k\},-k\}$, and $u_k$ solve the equation
\begin{equation*}
-\Delta u_k=f_k,\s u_k|_{\partial D}=0.
\end{equation*}
Then $\|f_k\|_{L^1\log L^1}\leq \|f\|_{L^1\log L^1}$,
$u_k\in W^{2,p}$ and
$$
\int_D|\nabla u_k|^2dx=\int_Dufdx
\leq \|u_k\|_{L^\infty(D)}\|f\|_{L^1(D)}
\leq \frac{1}{4}(\|u_k\|_{L^\infty(D)}+\|f_k\|_{L^1(D)})^2.
$$
Then we may assume $u_k$ converges weakly in $W^{1,2}$. Obviously,
the limit is $u$. 
It is easy to check that $\|f\|_{L^1}\leq C(1+\|f\|_{L^1\log L^1})$.
Thus we get \eqref{Hardy2.2}.
\endproof

It follows from Corollay \ref{Linfinity2}  the following:

\begin{cor}\label{oscandnablaL2} Let $q\in (0,2)$, $f\in L^1{\log L^1}(D)$
and $u$ solve the equation
\begin{equation*}\label{equationofu}
-\Delta u=f.
\end{equation*}
Then 
\begin{equation*}\label{Hardy3}
\|u-c\|_{L^\infty(D_\frac{1}{2})}+\|\nabla u\|_{L^2(D_\frac{1}{2})}
\leq C(1+\|f\|_{L^1{\log L^1}(D)}+\|\nabla u\|_{L^q(D)}),
\end{equation*}
where $c$ is the mean value of $u$ over $D_\frac{1}{2}$.
\end{cor}

\proof Let $\phi$ be a cut-off function, which is 1 on $D_\frac{1}{2}$
and $0$ on $D\setminus D_\frac{3}{4}$. Let $v$ solve the equation
$$
-\Delta v=2\nabla (u-c)\nabla \phi+(u-c)\Delta \phi,\s
v|_{\partial v}=0.
$$
By Poincare inequality and $L^p$-estimate,  $\|v\|_{W^{2,q}(D)}
<C\|\nabla u\|_{L^q(D)}$. Since $-\Delta(\phi(u-c)-v)=f$,
by Corollary \ref{Linfinity2}, we get the result.
\endproof

\section{$\epsilon$-regularity}

The basic tool of us is the Gauss equation on an isothermal
coordinate system. So, we need to discuss such type of equation first.
We will base on the following:
\begin{thm}\cite{Brezis-Merle}\label{Briezis} Let $u$ be a solution of 
$$
-\Delta u=f,\s u|_{\partial\Omega}=0,
$$
where $f\in L^1$ and $\Omega$ is a bounded domain of $\R^2$. Then for any $\epsilon>0$, we have
$$\int_\Omega e^{\frac{(4\pi-\epsilon)|u|}{\|f\|_{L^1(\Omega)}}}<\frac{4\pi
diam(\Omega)}{\epsilon^2}.$$\\
\end{thm}

As a corollary, we have
\begin{pro}\label{key} Let $u$ solve the equation $-\Delta u=Ke^{2u}$ on the unit 2-dimensional
disk $D$ with
$$\int_D|K|\log(1+|K|)e^{2u}dx<\Lambda.$$
Then  there exists an $\epsilon_0$, such that
if 
$$
\int_De^{2u}dx\leq\epsilon<\epsilon_0,
$$
then
$$
\int_{D_\frac{1}{2}}|K|e^{2u}\log(1+|K|e^{2u})dx\leq
C(\Lambda,\epsilon_0).
$$
\end{pro}

\proof

Let $v$ be the solution of
$$-\Delta v=Ke^{2u},\, v|_{\partial D}=0.$$
Applying Theorem \ref{Briezis} we have
$$\int_D e^{2s|v|}dx<C_1,$$
whenever $2s\leq 2s_0<\frac{4\pi}{\|Ke^{2u}\|_{L^1(D)}}$.

Since
\begin{eqnarray*}
\int_D|K|e^{2u}dx&<&
\int_{\{|K|>A\}}|K|e^{2u}dx+\int_{\{|K|\leq A\}}|K|e^{2u}dx\\
&\leq&\frac{1}{\log(1+A)}\int_D|K|e^{2u}\log(1+|K|)dx+
A\epsilon_0\\
&\leq&\frac{\Lambda}{\log(1+A)}+A\epsilon_0,
\end{eqnarray*}
we may choose $A$ to be sufficiently
large such that $\frac{\Lambda}{\log(1+A)}<\frac{\pi}{6}$,
and choose $\epsilon_0$ to be sufficiently small such that 
$A\epsilon_0<\frac{\pi}{6}$. Thus
$$\int_De^{6|v|}dx<C_2,$$
which implies that 
$$\int_D|v|<C_3.$$

By Jensen's inequality
$$
\int_{D_\frac{1}{4}(x)}2udx\leq \log\int_{D_\frac{1}{4}
(x)}e^{2u}dx\leq C_4,\s \forall x\in D_\frac{3}{4}.
$$
Since
$-\Delta(u-v)= 0$, for any $x\in D_\frac{3}{4}$, 
$$
u(x)-v(x)=\frac{1}{|D_\frac{1}{8}|}
\int_{D_\frac{1}{8}(x)}(u-v)dx\leq C_5=\frac{1}{4}C_4+C_3.
$$

Hence, we get
$$
u(x)\leq v(x)+C_5,\s \forall x\in D_\frac{3}{4}.
$$
Then
$$
\int_{D_\frac{3}{4}}e^{6u}\leq\int_{D_\frac{3}{4}}e^{6v+6C_5}dx\leq C_6,
$$
hence
\begin{eqnarray}\nonumber
\int_D|K|e^{2u}\log(1+|K|e^{2u})dx&<&
\int_D|K|e^{2u}\log(1+|K|)dx+\int_D|K|e^{2u}\log(1+e^{2u})dx\\ \nonumber
&<&\Lambda+(\int_{e^{2u}<|K|}+\int_{e^{2u}\geq |K|})|K|e^{2u}\log(1+e^{2u})dx\\\label{upper}
&\leq&2\Lambda+\int_{\{e^{2u}\geq |K|\}}e^{4u}\log(1+e^{2u})\\ \nonumber
&\leq&2\Lambda+\int_De^{6u}\\\nonumber
&<&C_7.
\end{eqnarray}

Set $v'\in W_0^{1,2}(D_\frac{3}{4})$, which solves the equation
$$-\Delta v'=Ke^{2u}.$$
By Corollary \ref{Linfinity2}, $\|v'\|_{C^0(D_\frac{3}{4})}<C_8$. Since we also have
$$
u(x)-v'(x)=\frac{1}{|D_\frac{1}{8}|}
\int_{D_\frac{1}{8}(x)}(u-v)d\sigma,
$$
we have
$$u(x)\leq C_9$$
for any $x\in D_\frac{1}{2}$. Thus we get from \eqref{upper} that
$$
\int_{D_\frac{1}{2}}|K|e^{2u}\log(1+|K|e^{2u})dx<2\Lambda+\int_{D_\frac{1}{2}}
e^{6(u-v)+6\|v\|_{L^\infty(D_\frac{1}{2})}} <C_{10}.
$$

\endproof

\begin{cor} Let $u$ solve the equation $-\Delta u=Ke^{2u}$ on the unit 2-dimensional
disk $D$ with
$$\int_D|K|^pe^{2u}dx<\Lambda.$$
Then  there exists an $\epsilon_0'$, such that
if 
$
\int_De^{2u}dx\leq\epsilon<\epsilon_0',
$
then
$$
\int_{D_\frac{1}{2}}|K|^pe^{2pu}dx\leq
C(p,\Lambda,\epsilon_0').
$$
\end{cor}

\proof
It is easy to check that $|K|\log(1+|K|)<C+|K|^p$. 
Set $v\in W_0^{1,2}(D_\frac{3}{4})$, which solves the equation
$$-\Delta v=Ke^{2u}.$$
By Proposition \ref{key} and Corollary \ref{Linfinity2},
we can choose $\epsilon_0'$ to be sufficiently small such that
$\|v\|_{L^\infty(D)}<C$.  Following the argument in the proof of 
Proposition \ref{key}, we get
$$u(x)\leq C$$
for any $x\in D_\frac{1}{2}$. Thus we get
$$
\int_{D_\frac{1}{2}}|K|^pe^{2pu}
\leq \int_{D_\frac{1}{2}}|K|^pe^{2u_k}d\sigma e^{2(p-1)\sup_{D_\frac{1}{2}}u}<C.
$$

\endproof

For  $g=e^{2u}g_{\S2}$, we define
$$
\rho(u,x)=\inf\{r:\int_{B_r^{g_{\S2}}(x)}e^{2u}d\mu_{g_{\S2}}=\frac{\epsilon_1}{2}\},
$$
and
$$
\rho(g)=\inf_{x\in\S2}\rho(u,x).
$$\\

\begin{cor}\label{convergence}
Let $g_k=e^{2u_k}g_{\S2}\in\mathcal{M}_\Lambda$.  We assume 
$\mu(g_k)>\tau$ for some $\tau >0$. If $\rho(g_k,x)>a>0$
for  a fixed sufficiently small $\epsilon_1$, then $u_k$
is bounded in $L^\infty$ and converges weakly in $W^{1,2}$.
\end{cor}

\proof For any $x_0\in\mathbb{S}^2$, we can choose conformal diffeomorphism $\phi_{x_0}:D\rightarrow B_{a}(x_0)$, such that $\phi_{x_0}(0)=x_0$.
Then $\phi_{x_0}$ defines an isothermal coordinate systerm with $x_0=0$.
It is follows from  Lemma \ref{key} that $\int_{\S2}|K_{g_k}e^{2u_k}|\log(1+|K_{g_k}|e^{2u_k})
d\mu_{g_{\S2}}<C$.  By Lemma \ref{Lq} and Corollary \ref{oscandnablaL2}, 
$\|\nabla u_k\|_{L^2}<C$.
By Poincare inequality and Sobolev embedding inequality,  
$\|u_k-\bar{u}_k\|_{L^1}<C$, where $\bar{u}_k$ is the mean value of
$u_k$.   Then applying Lemma \ref{Lq} and Corollary \ref{oscandnablaL2}
again,  we get 
$\|u_k-\bar{u}_k\|_{L^\infty}<C$. Since $\tau\leq \int e^{2u_k}=e^{2\bar{u}_k}\int e^{2(u_k-\bar{u}_k)}\leq
\Lambda$, we get $|\bar{u}_k|<C$.

\endproof

\section{A sequence $\{g_k\}$ with $\|K_g-1\|_{L^1}\rightarrow 0$}

The main task of this section is to prove the following:

\begin{lem}\label{sphere}
Let $g_k=e^{2u_k}g_{\S2}$. We assume $\mu(g_k)\leq\Lambda$ and
$\|K_{g_k}-1\|_{L^1(S^2,g_k)}\rightarrow 0$. After passing to a subsequence,
we can find a Mobious
transformation $\sigma_k$, such that $\sigma_k^*(g_k)$
converges to $g_{\S2}$ weakly in $W^{1,2}$. Moreover, we have
$\sigma_k^*(g_k)=e^{2u_k'}g_{\S2}$ with $|u_k'|<C$.
\end{lem}

\proof
Since 
\begin{eqnarray*}
\mu(g_k)&=&\int_{S^2}(1-K_{g_k})d\mu_{g_k}+\int_{S^2}K_{g_k}d\mu_{g_k}\\
&=& \int_{S^2}(1-K_{g_k})d\mu_{g_k}+4\pi,
\end{eqnarray*}
we get 
\begin{equation}\label{area}
\mu(g_k)\rightarrow 4\pi.
\end{equation}

We prove the lemma by contradiction.  By Corollary \ref{convergence},
we may suppose there exists $x_k$, such that 
$$
\rho(u_k,x_k)=\rho(g_k)\rightarrow 0.
$$
Let $y_k$ be the antipodal point of $x_k$, and $\pi_k$ be the 
stereprojection from $\S2\setminus\{y_k\}$ to $\C$. 
It is well-known that $\pi_k$
defines an isothermal coordinate system with $x_k=0$. 
Set 
$$
g_k=e^{2v_k}dz\otimes d\bar{z}.
$$
We have  $-\Delta v_k=K_{g_k}e^{2v_k}.$

Set $t_k$,
such that $\pi_k(B_{\rho(u_k,x_k)}(x_k))=D_{t_k}(0)$. Then
$t_k\rightarrow 0$.  
By the definition of $x_k$, 
we have
$$
\int_{D_{t_k}(0)}e^{2v_k}=\epsilon_1.
$$
Moreover, for any $R$, we can find $\tau(R)>0$, such that 
$\pi_k^{-1}(D_{\tau(R)t_k}(x))\subset B_{t_k}(\pi^{-1}(x))$ holds for
any $x\in D_R$. Thus,
$$
\int_{D_{\tau(R)t_k}(x)}e^{2v_k}
\leq \int_{B_{\rho(u_k,x_k)}(\pi_k^{-1}(x))}e^{2u_k}d\mu_{g_{\S2}}<\epsilon_1,\s \forall x\in D_R(0).
$$
By Lemma \ref{Lq}, we also have
$$
t^{2-q}\int_{D_t(x)}|\nabla v_k|^q<C(q,R),\s \forall x\in B_R(0).
$$

Let $v_k'(x)=v_k(t_kx)+\log t_k$.  Then, we have
$$
-\Delta v_k'=K_k(t_kx)e^{2v_k'},
$$
$$
\int_{D_R(0)}e^{2v_k'}=\int_{D_{Rt_k}(0)}e^{2v_k'}\leq \mu(g_k),\s \int_{D_1(0)}e^{2v_k'}=\int_{D_{t_k}(0)}e^{2v_k}= \epsilon_1,
$$
$$
\int_{D_{\tau(R)}(x)}e^{2v_k'}=\int_{D_{\tau(R)t_k}
(x_k+\tau(R)x)}e^{2v_k}<\epsilon_1, \forall x\in D_{R/t_k}(0),
$$
and
$$
\int_{D_R}|\nabla v_k'|^q=R^{2-q}(t_kR)^{q-2}\int_{D_{Rt_k}}|\nabla v_k|^q<C(R,q).
$$
Let $c_k$ be the mean value of $v_k'$ on $D_1$. 
By Poincare inequality and Corollary \ref{oscandnablaL2}, we get
$\|v_k-c_k\|_{L^\infty(B_R)}<C(R)$. Since $\epsilon_1=
\int_{D_1}e^{2v_k'}$, $|c_k|$ is bounded. Therefore, we may
assume $v_k'$ converges to a function $v$ weakly in $W^{1,2}$ and
$\|v_k'\|_{L^\infty(D_R)}<C(R)$.
Since
$$
\int_{D_R}|K_k-1|e^{2v_k'}\leq \|K_k-1\|_{L^1(S^2,g_k)}\rightarrow 0,
$$ 
$v$ satisfies the equation
$$
-\Delta v=e^{2v},\s \int_{\R^2}e^{2v}\leq 4\pi.
$$
By Lemma \ref{Ding} or Theorem \ref{Chen-Li},  
$$
e^{2v}dz\otimes d\bar{z}=g_{\S2},\s and \s \int_{\R^2}e^{2v}=4\pi.
$$ 
We usally call $v$ a bubble of $g_k$. 
 
Without loss of generality, we may assume $v_k'$
converges to $v$ almost everhwhere. Since $e^{v_k'}$ is bounded in
$L^p$ for any $p>2$, we get
$$
\int_{D_R}e^{2v_k'}dx\rightarrow\int_{D_R}e^{2v}dx.
$$

Let $x_0$ be the limit of $x_k$.
It is easy to check that
$$
\mu(g_k,\S2\setminus B_\delta(x_0))\geq\int_{D_{Rt_k}}e^{2v_k}=
\int_{D_R}e^{2v_k'}\rightarrow\int_{D_R}e^{2v},
$$
therefore,
$$
\lim_{k\rightarrow+\infty}\mu(g_k)\geq
\lim_{\delta\rightarrow 0}\varliminf_{k\rightarrow+\infty}
\mu(g_k,\S2\setminus B_\delta(x_0))+4\pi.
$$

Note that $\phi_k(x)=\pi_k^{-1}(\pi_k(x)/t_k)$ can be extended to a Mobious 
transformation $\sigma_k$ defined on $\S2$. Set
 $\sigma_k^*(g_k)=e^{2u_k'}g_{\S2}$. Then $u_k'$ converges 
weakly in $W^{1,2}(\S2\setminus\{y_0\})$, 
and $\|u_k'\|_{L^\infty(\S2\setminus B_\delta(y_0))}<C(\delta)$,
where $y_0$ is the  antipodal point of $x_0$. 

In fact, $u_k'$ must converges weakly in $W^{1,2}(\S2)$, and be bounded
in $L^\infty(\S2)$. 
Otherwise,  we can also get a bubble of $g_k'$ and imply that 
$$
\lim_{k\rightarrow+\infty}\mu(g_k')\geq
\lim_{\delta\rightarrow 0}\varliminf_{k\rightarrow+\infty}
\mu(g_k',\S2\setminus B_\delta(y_0))+4\pi
=\mu(g_{\S2})+4\pi= 8\pi.
$$
This  contradicts \eqref{area}.
\endproof

\section{The proof of the Theorem \ref{main}}
It is easy to deduced Theorem \ref{main} from Corollary \ref{convergence} and the following lemma:

\begin{lem}
For any
$\Lambda$, there exits $\tau$ and $a>0$,  such that if $g\in\mathcal{M}_\Lambda$ and 
$\|K_g-1\|_{L^1(S^2,g)}<\tau$, then we can find a Mobius transformation $\sigma$, such that 
$\rho(\sigma^*(g))>a$.
\end{lem}

\proof
Suppose the lemma is not true. Then we can find $\lambda_k\rightarrow 0$,
$g_k\in\mathcal{M}_\Lambda$, such that
for any $\sigma_k$, $\rho(\sigma_k^*(g_k))\rightarrow 0$.

However, by Lemma \ref{sphere}, we can find
$\sigma_k$,
such that $\sigma_k^*(g_k)=e^{2u_k}$ with $\|u_k\|_{L^\infty}<C$, and
$u_k$ converges weakly in $W^{1,2}$, thus
$\rho(\sigma_k^*(g_k))$ can not converges to 0. 
\endproof

\end{document}